\documentclass[fleqn,12pt,twoside]{article}
\usepackage{espcrc1}
\usepackage{multirow}
\usepackage{array,booktabs}

\usepackage{graphicx}
\usepackage[figuresright]{rotating}

\newcommand{\AmS}{{\protect\the\textfont2
  A\kern-.1667em\lower.5ex\hbox{M}\kern-.125emS}}

\title{An effective numerical method to solve a class of nonlinear singular boundary value problems using improved differential transform method}

\author{Lie-jun Xie \thanks{Corresponding author. } \thanks{E-mail address: xieliejun@nbu.edu.cn, xieliejun@hotmail.com.}, Cai-lian Zhou and Song Xu \\
        \vskip 10pt {Department of Mathematics, Faculty of Science,
        Ningbo University \\
        Ningbo Zhejiang, 315211 P. R. China}
       }

\begin{document}

\maketitle

\begin{abstract}
In this work, an effective numerical method is developed to solve a class of singular boundary value problems arising in various physical models by using the improved differential transform method (IDTM). The IDTM applies the Adomian polynomials to handle the differential transforms of the nonlinearities arising in the given differential equation. The relation between the Adomian polynomials of those nonlinear functions and the coefficients of unknown truncated series solution is given by a simple formula, through which one can easily deduce the approximate solution which takes the form of a convergent series. An upper bound for the estimation of approximate error is presented. Several physical problems are discussed as illustrative examples to testify the validity and applicability of the proposed method. Comparisons are made between the present method and the other existing methods.


\end{abstract}

\section{Introduction}
Singular boundary value problems (SBVPs) is an important class of boundary value problems, and arises frequently in the modeling of many actual problems related to physics and engineering areas such as in the study of electro hydrodynamics, theory of thermal explosions, boundary layer theory, the study of astrophysics, three layer beam, electromagnetic waves or gravity driven flows, inelastic flows, the theory of elastic stability and so on. In general, SBVPs is difficult to solve analytically. Therefore, various numerical techniques have been proposed to treat it by many researchers. However, the solution of SBVPs is numerically challenging due to the singularity behavior at the origin.

In this work, we are interested again in the following SBVPs arising frequently in applied science and engineering:
\begin{equation} \label{1-1}
\hskip 25pt
\displaystyle u''(x)+\frac{\alpha}{x}u'(x)=f(x,u), \hskip 3pt 0<x\leq 1, \hskip 3pt \alpha \geq 1,
\end{equation}
subject to the boundary value conditions
\begin{equation}  \label{1-2}
\hskip 25pt
u'(0)=0
\end{equation}
and
\begin{equation}  \label{1-3}
\hskip 25pt au(1)+bu'(1)=c,
\end{equation}
where $a,b$ and $c$ are any finite real constants. If $\alpha=1$, (1) becomes a cylindrical problem, and it becomes a spherical problem when $\alpha=2$. It is assumed that $f(x,u)$ is continuous, $\frac{\partial f}{\partial u}$ exists and is continuous and $\frac{\partial f}{\partial u}\geq 0$ for any $0<x\leq 1$ such that equation (1) has a unique solution \cite{Russel75}.

The SBVPs (1)-(3) with different $\alpha$ arise in the study of various scientific problems for certain linear or nonlinear functions $f(x,u)$. The common cases related to the actual problems are summarized as follows. The first case for $\alpha=2$ and
\begin{equation}
\hskip 25pt f(x,u)=f(u)=\frac{\delta u(x)}{u(x)+\mu}
\end{equation}
emerges from the modeling of steady state oxygen diffusion in a spherical cell with Michaelis-Menten uptake kinetics \cite{Lin76,McElwain78}. In this case, $u(x)$ represents the oxygen tension; $\delta$ and $\mu$ are positive constants involving the reaction rate and the Michaelis constant. Hiltmann and Lory \cite{Hiltmann83} proposed the existence and uniqueness of the solution for $b=1$ and $a=c$. Analytical bounding functions were given in \cite{Anderson85}. The numerical methods to solve the SBVPs for this case have attracted a reasonable amount of research works, such as the finite difference method (FDM) \cite{Pandey97}, the cubic spline method (CSM) \cite{Rashidinia07,Kanth06}, the Sinc-Galerkin method (SGM) \cite{Babolian13}, the Adomian decomposition method (ADM) and its modified methods \cite{Khuri10,Wazwaz13,Singh14}, the variational iteration method (VIM) \cite{Kanth10,Wazwaz11}, the series expansion technique (SEM) \cite{Turkyilmazoglu13} and the B-spline method (BSM) \cite{Caglar09}.

The second case arises in the study of the distribution of heat sources in the human head \cite{Flesch75,Gray80,Duggan86} , in which $\alpha=2$ and
\begin{equation}
\hskip 25pt f(x,u)=f(u)=-le^{-lku(x)}, l>0, k>0.
\end{equation}
In \cite{Duggan86}, point-wise bounds and uniqueness results were presented for the SBVPs with the nonlinear function $f(x,u)$ of the forms given by (4) and (5). Quite a little amount of works by using different approaches, including the FDM \cite{Pandey97}, the CSM \cite{Rashidinia07,Kanth06} and the SGM \cite{Babolian13}, have been proposed to obtain the approximate solutions of this case.

The third important case of physical significance is when $\alpha=1, 2$ and
\begin{equation}
\hskip 25pt f(x,u)=f(u)=\nu e^{u(x)},
\end{equation}
which arises in studying the theory of thermal explosions \cite{Khuri10,Kumar10,Chang14} and the electric double layer in a salt-free solution \cite{Chang12}. A variety of numerical methods have been applied to handle such SBVPs, for example, the fourth order finite difference method (FFDM) \cite{Chawla88}, the modified Adomian decomposition method \cite{Khuri10,Singh14,Kumar10}, the Taylor series method (TSM) \cite{Chang14} and the BSM \cite{Caglar09}.

Besides, Chandrasekhar \cite{Chandrasekhar39} derived another case for $\alpha=2, b=0$ and
\begin{equation}
\hskip 25pt f(x,u)=f(u)=-u^{\gamma}(x),
\end{equation}
which $\gamma$ is a physical constant. This case is in connection with the equilibrium of thermal gas thermal \cite{Ames68}. The numerical solution of this kind of equation for  $\gamma=5$ was considered by using various methods, such as the FFDM \cite{Chawla88}, the VIM \cite{Kanth10}, the SEM \cite{Turkyilmazoglu13} and the modified Adomian decomposition method \cite{Singh14}.

All the aforementioned methods can yield a satisfied result. However, each of these methods has its own weaknesses. For example, the VIM \cite{Kanth10,Wazwaz11} has an inherent inaccuracy in identifying the Lagrange multiplier, and fails to solve the equation when the nonlinear function $f(x,u)$ is of the forms (5) and (6). Those methods such as the FDM \cite{Pandey97,Chawla88}, the SEM \cite{Turkyilmazoglu13}, the SGM \cite{Babolian13} and the spline method \cite{Rashidinia07,Kanth06,Caglar09} require a tedious process and huge volume of computations in dealing with the linearization or discretization of variables. The ADM \cite{Wazwaz13} needs to obtain the corresponding Volterra integral form of the given equation, via which one can overcome the difficulty of singular behavior at $x=0$. The modified ADM \cite{Khuri10,Kumar10} needs to introduce a twofold indefinite integral operator to give better and accurate results; moreover, the success of method in \cite{Singh14} relies on constructing Green's function before establishing the recursive relation for applying the ADM to derive the solution components. All those manners are at the expense of computation budgets. Besides, none of above methods is applied to handle the equations with all forms of nonlinearities (4)-(7).

In recent years, a lot of attentions have been devoted to the applications of differential transform method (DTM) and its modifications. The DTM proposed by Pukhov \cite{Pukhov80,Pukhov82,Pukhov86} at the beginning of 1980s. However, his work passed unnoticed. In 1986, Zhou reintroduced the DTM to solve the linear and nonlinear equations in electrical circuit problems \cite {Zhou86}. The DTM is a semi-numerical-analytic method that generates a Taylor series solution in the different manner. In the past forty years, the DTM has been successfully applied to solve a wide variety of functional equations; see \cite{Xie15} and the references therein.

Although being powerful, there still exist some difficulties in solving various of equations by the classical DTM. Some researchers have devoted to deal with these obstacles so as to extend the applications of the DTM. For example, in view of the DTM numerical solution cannot exhibit the real behaviors of the problem, Odibat {\it et al}. \cite{Odibat10} proposed a multi-step DTM to accelerate the convergence of the series solution over a large region and applied successfully to handle the Lotka-Volterra, Chen and Lorenz systems. In \cite{Gokdogan12,Momani08}, authors suggested an alternative scheme to overcome the difficulty of capturing the periodic behavior of the solution by combining the DTM, Laplace transform and Pad\'{e} approximants. Another difficulty  is to compute the differential transforms of the nonlinear components in a simple and effective way. By using the traditional approach of the DTM, the computational difficulties will inevitably arise in determining the transformed function of an infinity series. Compared to the traditional method, Chang and Chang \cite{Chang08} proposed a relatively effective algorithm for calculating the differential transform through a derived recursive relation. Yet, by using their method, it is inevitable to increase the computational budget, especially in dealing with those differential equations which have two or more nonlinear terms being investigated. Recently, the authors \cite{Ebaid12,Fat13} disclosed the relation between the Adomian polynomials and the differential transform of nonlinearities, and developed an inspiring approach to handle the nonlinear functions in the given functional equation. Meanwhile, the problem of tedious calculations in dealing with nonlinear problems by using the ADM has also been improved considerably by Duan \cite{Duan10-6,Duan10-7,Duan11}. All of these effective works make it possible to broaden the applicability and popularity of the DTM considerably.

The aim of this work is to develop an efficient approach to solve the SBVPs (1)-(3) with those nonlinear terms (4)-(7). This scheme is mainly based on the improved differential transform method (IDTM), which is the improved version of the classical DTM by using the Adomian polynomials to handle the differential transforms of those nonlinear functions (4)-(7). No specific technique is required in dealing with the singular behavior at the origin. Meanwhile, unlike some existing approaches, the proposed method tackles the problem in a straightforward manner without any discretization, linearization or perturbation. The numerical solution obtained by the proposed method takes the form of a convergent series with those easily computable coefficients through the Adomian polynomials of those nonlinear functions as the forms of (4)-(7).

The rest of the paper is organized as follows. In the next section, the concepts of DTM and Adomian polynomials are introduced. Algorithm for solving the problem (1)-(3) and an upper bound for the estimation of approximate error are presented in Section 3. Section 4 shows some numerical examples to testify the validity and applicability of the proposed method. In Section 5, we end this paper with a brief conclusion.

\section{Adomian polynomial and Differential transform}

\subsection{Adomian polynomial}
In the Adomian decomposition method (ADM), a key notion is the Adomian polynomials, which are tailored to the particular nonlinearity to easily and systematically solve nonlinear differential equations. The interested readers are referred to Refs. \cite{Adomian90,Adomian94} for the details of the ADM.

For the applications of decomposition method, the solution of the given equation in a series form is usually expressed by
\begin{equation}
\hskip 25pt u= \sum_{m=0}^{\infty} u_m,
\end{equation}
and the infinite series of polynomials
\begin{equation}
\hskip 25pt f(u)=f \bigg{(} \sum_{m=0}^{\infty} u_m \bigg{)} =\sum _{m=0}^{\infty} A_{m}
\end{equation}
for the nonlinear term $f(u)$, where $A_{m}$ is called the Adomian polynomials, and depends on the solution components $u_0, u_1, \cdots, u_m$. The traditional algorithm for evaluating the Adomoan polynomials $A_n$ was first provided in \cite{Adomian83} by the formula
\begin{equation}
\hskip 25pt A_{n}= \frac {1}{n!} \frac {\mathrm{d}^n}{\mathrm{d} \lambda^n} f \bigg{(}\sum_{m=0}^{\infty} u_m \lambda^m \bigg{)} \bigg{|}_{\lambda=0}.
\end{equation}
A large amount of works \cite{Duan10-6,Duan10-7,Duan11,Adomian83,Rach08,Rach84,Wazwaz00,Abbaoui95,Abdelwahid03,Azreg09} have been applied to give the more effective computational method for the Adomian polynomials. For fast computer generation, we favor Duan's Corollary 3 algorithm \cite{Duan11} among all of these methods, as it merely involves the analytic operations of addition and multiplication without the differentiation operator, which is eminently convenient for symbolic implementation by computer algebraic systems such as Maple and Mathematics. The method to generate the Adomian polynomials in \cite{Duan11} is described as follows:
\begin{equation}
\begin{array}{ll}
\hskip 25pt C_{n}^1= u_n, \hskip 10pt n \geq 1,\\
\hskip 25pt \displaystyle C_{n}^k= \frac{1}{n} \sum_{j=0}^{n-k} (j+1)u_{j+1}C_{n-1-j}^{k-1}, \hskip 10pt  2\leq k \leq n,
\end{array}
\end{equation}
such that
\begin{equation}
\begin{array}{ll}
\hskip 25pt A_0= f(u_0),\\
\hskip 25pt \displaystyle A_n= \sum_{k=1}^{n} C_{n}^{k} f^{(k)}(u_0), \hskip 10pt  n \geq 1.
\end{array}
\end{equation}
It is worth mentioning that Duan's algorithm involving (11) and (12) has been testified to be one of the fastest subroutines on record \cite{Duan11}, including the fast generation method given by Adoamin and Rach \cite{Adomian83}.

\subsection{Differential transform}
The differential transform of the $k^{th}$ differentiable function $u(x)$ at $x=0$ is defined by
\begin{equation}
\hskip 25pt U(k)=\frac{1}{k!}\bigg{[} \frac{\mathrm{d}^k u(x)}{\mathrm{d}x^k}  \bigg{]}_{x=0},
\end{equation}
and the differential inverse transform of $U(k)$ is described as
\begin{equation}
\hskip 25pt u(x)=\sum_{k=0}^{\infty}U(k) x^k,
\end{equation}
where $u(x)$ is the original function and $U(k)$ is the transformed function.

For the practical applications, the function $u(x)$ is expressed by a truncated series and Eq. (14) can be written as
\begin{equation}
\hskip 25pt u(x)\approx u_N(x)=\sum_{k=0}^{N}U(k) x^k.
\end{equation}

It is not difficult to deduce the transformed functions of the fundamental operations listed in Table 1.
\begin{table}[htb]
\caption{The fundamental operations of the DTM.}
\label{table:1}
\newcommand{\m}{\hphantom{$-$}}
\newcommand{\cc}[1]{\multicolumn{1}{c}{#1}}
\renewcommand{\tabcolsep}{4pc} 
\renewcommand{\arraystretch}{1.2} 
\begin{tabular}{@{}ll}
\toprule
\small Original function      &  \small Transformed function \\
\hline
\small $w(x)=\alpha u(x)\pm \beta v(x)$       & \small $W(k)=\alpha U(k)\pm \beta V(k)$ \\
\small $w(x)=u(x)v(x)$       & \small $W(k)=\sum _{m=0}^k U(m)V(k-m)$ \\
\small $w(x)=\mathrm{d}^m u(x)/\mathrm{d}x^m$       & \small $W(k)= \frac{(k+m)!}{k!}U(k+m)$ \\
\small $w(x)=x^m$       & \small $W(k)= \delta(k-m)=\left \{
                                                        \begin{array}{ll}
                                                            1, \hskip 3pt \mathrm{if} \hskip 3pt k=m,\\
                                                            0, \hskip 3pt \mathrm{if} \hskip 3pt k\neq m.
                                                        \end{array}
                                                     \right.
                                                     $ \\
\small $w(x)=\exp(x)$       & \small $W(k)= 1/k!$ \\
\small $w(x)=\sin(\alpha x+\beta)$       & \small $W(k)= \alpha ^k/k! \sin(k \pi /2+\beta)$ \\
\small $w(x)=\cos(\alpha x+\beta)$       & \small $W(k)= \alpha ^k/k! \cos(k \pi /2+\beta)$ \\
\bottomrule
\end{tabular}\\[2pt]
Note that $\alpha, \beta$ are constants and $m$ is a nonnegative integer.
\end{table}

\section{Method of solution of SBVPs (1)-(3)}
We want to find the approximate solution of the problem (1)-(3) with the type:
\begin{equation}
\hskip 25pt u_N(x)=\sum_{k=0}^{N}U(k) x^k,
\end{equation}
where the coefficients $U(0),U(1),\cdots,U(N)$ are determined using the following steps:

$\bullet $ According to the definition (13) of the differential transform and the boundary value condition (2), we have
\begin{equation}
\hskip 25pt U(1)=0.
\end{equation}
Suppose that
\begin{equation}
\hskip 25pt U(0)=\beta,
\end{equation}
where $\beta$ is a real parameter to be determined.

$\bullet $ Multiplying both sides of Eq. (1) by variable $x$, we have
\begin{equation}
\hskip 25pt
\displaystyle xu''(x)+\alpha u'(x)=xf(x,u).
\end{equation}
Applying the differential transform (13) to Eq. (19), we get the following recurrence relation:
\begin{equation}
\hskip 25pt U(k+1)=\frac{F(k-1)}{(k+1)(k+\alpha)},  \hskip 5pt k=1,2,\cdots,N-1,
\end{equation}
where $F(k)$ is the differential transform of the nonlinear function $f(x,u)=f(u)$.

$\bullet $ Using Lemma 3.1 in \cite{Fat13}, we compute $F(k)$ through the Adomian polynomials $A_k$:
\begin{equation}
\hskip 25pt F(k)=A_k, \hskip 5pt k=0,1,2,\cdots,N.
\end{equation}
{\bf \noindent Remark 1.}  Lemma 3.1 in \cite{Fat13} indicates that the differential transforms and the Adomian polynomials of nonlinear functions have the same mathematical structure such that we can derive the differential transforms of any nonlinear functions by merely calculating the relevant Adomian polynomials but with constants instead of variable components.
{\bf \noindent Remark 2.} As mentioned before, we use Duan's Corollary 3 algorithm \cite{Duan11} (11)-(12) to generate the Adomian polynomials.

$\bullet $ Substituting (21) into (20), and then combining the relations (16)-(18), we obtain the truncated series solution of the problem (1)-(3) as follows:
\begin{equation}
\hskip 25pt u_N(x)=\beta+\sum_{k=1}^{N-1}\frac{A_{k-1}}{(k+1)(k+\alpha)} x^{k+1}.
\end{equation}

$\bullet $ Imposing the truncated series solution (22) on the boundary condition (\ref{1-3}), we obtain a nonlinear algebraic equation with unknown parameter $\beta$:
\begin{equation}
\hskip 25pt g(\beta)=0.
\end{equation}
Solving Eq. (23), and substituting the value of $\beta$ into (22), we obtain the final result.

An upper bound for the estimation of approximate error is presented in the following lemma.
\vskip 10pt
{\bf \noindent Lemma 1.} {\it
Suppose that $u(x) \in C^{N+1}[0,1]$ is the exact solution of the problem (1)-(3), $u_N(x)=\sum_{k=0}^{N}U(k)x^k$ is the truncated series solution with degree $N$, it holds that
\begin{equation}  \label{theo-err}
\hskip 25pt ||u(x)-u_N(x)||_{\infty} \leq \frac{M}{(N+1)!}+\max \limits_{0\leq k\leq N} \left | c_k  \right |,
\end{equation}
where $M=\max \limits_{0\leq x\leq1} | u^{(N+1)}(x)|$, $c_k= \frac{u^{(k)}(0)}{k!}-U(k)$. }
\vskip 10pt

{\bf \noindent Proof.}
Obviously, we have
\begin{equation}
\hskip 25pt ||u(x)-u_N(x)||_{\infty} \leq ||u(x)- \widetilde{u}_N(x)||_{\infty} +||\widetilde{u}_N(x)-u_N(x)||_{\infty} ,
\end{equation}
where $\widetilde{u}_N(x)=\sum_{k=0}^{N} \frac{u^{(k)}(0)}{k!}x^k$ is the Taylor polynomial of the unknown function $u(x)$ at $x=0$.

Since $u(x) \in C^{N+1}[0,1]$, it follows that
\[
\hskip 25pt u(x)=\widetilde{u}_N(x) + R_N(x)=\widetilde{u}_N(x) + \frac{u^{(N+1)}(\xi)}{(N+1)!}x^{N+1}, \hskip 5pt \xi \in (0,1),
\]
where $R_N(x)$ is the remainder of Taylor polynomial $\widetilde{u}_N(x)$.
Therefore
\begin{equation}
\hskip 25pt \left |  u(x)-\widetilde{u}_N(x)\right |=\left | R_N(x) \right |=\left | \frac{u^{(N+1)}(\xi)}{(N+1)!}x^{N+1} \right |
\leq \frac{1}{(N+1)!} \max \limits_{0\leq x\leq1} \left | u^{(N+1)}(x) \right |.
\end{equation}

Let
\[
\hskip 25pt {\bf C}= (c_0,c_1,\cdots,c_N), \hskip 5pt {\bf \Theta}=(x^0,x^1,\cdots,x^N)^T,
\]
where
\[
\hskip 25pt c_k= \frac{u^{(k)}(0)}{k!}-U(k), \hskip 5pt k=0,1,\cdots,N.
\]
We then have
\begin{equation}
\hskip 25pt \left |  \widetilde{u}_N(x)- u_N(x) \right |= \left | \sum_{k=0}^{N} \Bigg{(} \frac{u^{(k)}(0)}{k!}-U(k) \Bigg{)} x^k \right |
=\left |  {\bf C} \cdot {\bf \Theta} \right | \leq ||{\bf C}||_{\infty} \cdot ||{\bf \Theta}||_{\infty}
\end{equation}
Combining the relations (25)-(27), it follows that
\begin{equation}
\begin{array}{ll}
\hskip 25pt \displaystyle ||u(x)-u_N(x)||_{\infty} \leq  \frac{1}{(N+1)!} \max \limits_{0\leq x\leq1} \left | u^{(N+1)}(x) \right | + ||{\bf C}||_{\infty} \cdot ||{\bf \Theta}||_{\infty}\\
\hskip 120pt \displaystyle  \leq \frac{M}{(N+1)!}+ \max \limits_{0\leq k\leq N} \left | c_k  \right |.
\end{array}
\end{equation}
Thus, the proof is completed.

\section{Numerical examples}
In this section, based on the discussion in Section 3, we report numerical tests of five classical examples discussed frequently to testify the validity and applicability of the proposed method. All the numerical computations were performed using Maple and Matlab on personal computer. For comparison, we computed the absolute error defined by
\begin{equation}  \label {absol-err}
\hskip 25pt \mathrm{E}_N(x)=\left | u(x)-u_N(x)\right |
\end{equation}
and the maximal absolute error by
\begin{equation}  \label {max-absol-err}
\hskip 25pt \mathrm{ME}_N=\max \limits_{0 \leq x \leq 1} \left | u(x)-u_N(x)\right |,
\end{equation}
where $u(x)$ is the exact solution and $u_N(x)$ is the truncated series solution with degree $N$.

\vskip 10pt
{\bf \noindent Example 1.} Consider the following nonlinear SBVP in the study of isothermal gas sphere \cite{Singh14,Kanth10,Chawla88}:
\begin{equation}
\hskip 15pt
\displaystyle u''(x)+\frac{2}{x} u'(x)=-u^5(x),
\end{equation}
subject to the boundary conditions
\begin{equation}
\hskip 15pt  u'(0)=0, \hskip 5pt u(1)=\frac{\sqrt 3}{2}.
\end{equation}
The exact solution of this problem is given by $u(x)=\sqrt {\frac {3}{3+x^2}}$. It is also known as the Emden-Fowler equation of the first kind. In what follows, we shall solve it with the proposed algorithm.

Firstly, we set
\[
\hskip 25pt U(0)=\beta, \hskip 4pt U(1)=0.
\]
The Adomian polynomials of nonlinear term $f(x,u)=-u^5(x)$ in this problem are computed as

\[
\begin{array}{lll}
\hskip 25pt \displaystyle A_0=-U^5(0),\\
\hskip 25pt \displaystyle A_1=-5U^4(0)U(1) ,\\
\hskip 25pt \displaystyle A_2=-10U^3(0)U^2(1)-5U^4(0)U(2),\\
\hskip 25pt \displaystyle A_3=-10U^2(0)U^3(1)-20U^3(0)U(1)U(2)-5U^4(0)U(3),\\
\hskip 45pt \vdots
\end{array}
\]
Furthermore, according to the relations (20) and (21), we obtain the differential transforms $U(k)$ of the unknown function $u(x)$
\[
\begin{array}{lll}
\hskip 25pt \displaystyle U(2)=\frac{1}{2\cdot 3} A_0=-\frac{1}{6}\beta^5, \\
\hskip 25pt \displaystyle U(4)=\frac{1}{4\cdot 5} A_3=\frac{1}{24}\beta^9,\\
\hskip 25pt \displaystyle U(6)=\frac{1}{6\cdot 7} A_5=-\frac{5}{432}\beta^{13},\\
\hskip 45pt \vdots \\
\hskip 25pt U(k)=0, \hskip 5pt \mathrm{if} \hskip 5pt k \hskip 5pt \mathrm{is \hskip 5pt odd \hskip 5pt and } \hskip 5pt k \geq 3.
\end{array}
\]
By using Eq. (22), we obtain the truncated series solution for $N=10$ as follows:
\begin{equation}  \label{example1-1}
\hskip 25pt \displaystyle u_{10}(x)=\beta-\frac{1}{6}\beta^5x^2+\frac{1}{24}\beta^9x^4-\frac{5}{432}\beta^{13}x^6+\frac{35}{10368}\beta^{17}x^8-\frac{7}{6912}\beta^{21}x^{10}. \\
\end{equation}

Secondly, imposing the truncated series solution (\ref{example1-1}) on the boundary conditions $u(1)=\sqrt3/2$, we get a nonlinear algebraic equation. By solving it, the unknown parameter $\beta$ is computed as
\begin{equation}  \label{example1-2}
\hskip 25pt \beta=1.000553890.
\end{equation}

Finally, substituting (\ref{example1-2}) into (\ref{example1-1}), we get the approximate solution with degree $10$
\[
\begin{array}{lll}
\hskip 25pt \displaystyle u_{10}(x)=1.000553890-0.1671287533x^2+(0.4187483621e-1)x^4-\\
\hskip 72pt (0.1165769154e-1)x^6+(0.3407699551e-2)x^8-\\
\hskip 72pt (0.1024576736e-2)x^{10}.
\end{array}
\]

In Table 2, we compare the absolute errors (\ref {absol-err}) of numerical results obtained by the present method, the VIM \cite{Kanth10} and the modified ADM using Green functions (GIDM) \cite{Singh14} for $N=12$. Table 3 lists the theoretical estimate errors (\ref{theo-err}) and the maximal absolute errors (\ref {max-absol-err}) of the approximate solutions for changing approximation levels, and shows a comparison of the maximal absolute errors with the GIDM \cite{Singh14} and the FFDM \cite{Chawla88}. We can see from Table 3 that the accuracy of our computational results is getting better as the approximation level is increasing. Moreover, our numerical solution $u_{10}(x)$ has an accuracy of O($10^{-4}$), whereas the GIDM \cite{Singh14} needs to employ 14 terms to archive this goal as shown in table 1 of ref. \cite{Singh14}; numerical solution with even 64 terms obtained by the FFDM \cite{Chawla88} still hovers at this level. In summary, Tables 2 and 3 indicate that the results of our proposed method have higher accuracy than of the GIDM \cite{Singh14}, the FFDM \cite{Chawla88} and the VIM \cite{Kanth10}.

\begin{table}[htb]
\caption{Comparison of the absolute error $\mathrm{E}_{12}(x)$ for Example 1.}
\label{table:1}
\newcommand{\m}{\hphantom{$-$}}
\newcommand{\cc}[1]{\multicolumn{21}{c}{#1}}
\renewcommand{\tabcolsep}{3.1pc} 
\renewcommand{\arraystretch}{1.2} 
\begin{tabular}{@{}llll}
\toprule
$x$        & \small GIDM \cite{Singh14} & \small VIM \cite{Kanth10}  & \small Present method  \\
\hline
\small $0.0$        & \small 3.1880e-03         & \small 6.3220e-03    & \small 1.6776e-04  \\
\small $0.1$        & \small 3.1209e-03         & \small 6.2702e-03    & \small 1.6637e-04  \\
\small $0.2$        & \small 2.9269e-03         & \small 6.1173e-03    & \small 1.6227e-04   \\
\small $0.3$        & \small 2.6263e-03         & \small 5.8687e-03    & \small 1.5568e-04   \\
\small $0.4$        & \small 2.2489e-03         & \small 5.5281e-03    & \small 1.4691e-04    \\
\small $0.5$        & \small 1.8284e-03         & \small 5.0903e-03    & \small 1.3639e-04   \\
\small $0.6$        & \small 1.3978e-03         & \small 4.5347e-03    & \small 1.2450e-04  \\
\small $0.7$        & \small 9.8413e-04         & \small 3.8201e-03    & \small 1.1132e-04   \\
\small $0.8$        & \small 6.0707e-04         & \small 2.8837e-03    & \small 9.5269e-05  \\
\small $0.9$        & \small 2.7774e-04         & \small 1.6426e-03    & \small 6.8180e-05   \\
\small $1.0$        & \small 3.52e-08           & \small 1.00e-10      & \small 0          \\
\bottomrule
\end{tabular}\\[2pt]
\end{table}

\begin{table}[htb]
\caption{The theoretical estimate errors  $\mathrm{TE}_N$ and comparison of the maximal absolute errors $\mathrm{ME}_N$ of present method and of other methods for Example 1.}
\label{table:1}
\newcommand{\m}{\hphantom{$-$}}
\newcommand{\cc}[1]{\multicolumn{21}{c}{#1}}
\renewcommand{\tabcolsep}{0.54pc} 
\renewcommand{\arraystretch}{1.2} 
\begin{tabular}{@{}llllllllll}
\toprule
$N$      & \small $\mathrm{TE}_N$ & \small $\mathrm{ME}_N$ & \small $N$  & \small $\mathrm{TE}_N$ & \small $\mathrm{ME}_N$  & \small $N$   & \small in \cite{Singh14} & \small $N$   & \small in \cite{Chawla88}\\
\hline
\small $6$  & \small 1.83e-02 & \small 6.80e-03  & \small $12$ & \small 4.7721e-04 & \small 1.6776e-04  & \small $12$ & \small 1.3978e-03 & \small $16$ & \small 3.64e-04 \\

\small $8$  & \small 5.10e-03 & \small 1.70e-03  & \small $16$ & \small 4.6453e-05 & \small 1.6521e-05   & \small $16$ & \small 2.4654e-04 & \small $32$ & \small 2.49e-04   \\

\small $10$ & \small 1.5666e-03 & \small 5.5389e-04 & \small $20$ & \small 4.6453e-06 & \small 1.6614e-06  & \small $20$  & \small 4.8643e-05 & \small $64$ & \small 1.60e-04  \\
\bottomrule
\end{tabular}\\[2pt]
\end{table}

\vskip 10pt
{\bf \noindent Example 2.} Consider the following nonlinear SBVP \cite{Khuri10,Singh14,Caglar09,Chawla88}:
\begin{equation}
\hskip 15pt
\displaystyle u''(x)+\frac{1}{x} u'(x)=-e^{u(x)}, \\
\end{equation}
subject to the boundary conditions
\begin{equation}
\hskip 15pt
u'(0)=0, \hskip 5pt u(1)=0.
\end{equation}
The exact solution is given by $u(x)=2 \ln \frac{C+1}{Cx^2+1} $, where $C=3-2\sqrt2$.

The Adomian polynomials of nonlinear term $f(x,u)=-e^{u(x)}$ in this problem are computed as
\[
\begin{array}{lll}
\hskip 25pt \displaystyle A_0=-e^{U(0)},\\
\hskip 25pt \displaystyle A_1=-U(1)e^{U(0)},\\
\hskip 25pt \displaystyle A_2=-U(2)e^{U(0)}-\frac{1}{2}U^2(1)e^{U(0)},\\
\hskip 25pt \displaystyle A_3=-U(3)e^{U(0)}-U(1)U(2)e^{U(0)}-\frac{1}{6}U^3(1)e^{U(0)},\\
\hskip 45pt \vdots
\end{array}
\]

A comparison of the absolute errors (\ref{absol-err}) of the numerical solutions for $N=10, 20, 40$ obtained by the present method and the modified decomposition method (BSDM) \cite{Khuri10} is described in Table 4. Table 5 lists the maximal absolute errors (\ref{max-absol-err}) of those numerical results derived from the proposed method, the BSM \cite{Caglar09} and the FFDM \cite{Chawla88}. And also, we list the theoretical estimate errors (\ref{theo-err}) in Table 5 for comparison. It can be seen from Tables 4 and 5 that one can obtain the better approximate solution by using the present method compared to the other mentioned methods, even if we take the relative smaller $N$. Moreover, the theoretical estimate errors, the absolute errors and the maximal absolute errors all decrease as the increase of $N$. Therefore, evaluation of more components of the numerical solution will reasonably improve the accuracy.

\begin{table}[htb]
\caption{Comparison of the absolute errors $\mathrm{E}_N(x)$ for Example 2.}
\label{table:1}
\newcommand{\m}{\hphantom{$-$}}
\newcommand{\cc}[1]{\multicolumn{21}{c}{#1}}
\renewcommand{\tabcolsep}{1pc} 
\renewcommand{\arraystretch}{1.1} 
\begin{tabular}{@{}lllllll}
\toprule
    & \small {BSDM \cite{Khuri10} }       &        &         & \small {Present method}   &  &  \\
\cmidrule(r){2-4}  \cmidrule(r){5-7}

$x$                       & \small $\mathrm{E}_{10}(x)$         & \small $\mathrm{E}_{20}(x)$       & \small $\mathrm{E}_{40}(x)$        & \small $\mathrm{E}_{10}(x)$   & \small $\mathrm{E}_{20}(x)$  & \small $\mathrm{E}_{40}(x)$ \\
\hline

\small $0.0$              & \small  1.05e-05   & \small 1.05e-05  & \small  1.05e-05    & \small 1.05e-05  & \small 2.2e-09  & \small 1.4e-09 \\

\small $0.1$              & \small  1.05e-05   & \small 1.05e-05  & \small  1.05e-05    & \small 1.05e-05  & \small 1.2e-09  & \small 4.0e-10 \\

\small $0.2$              & \small  1.03e-05   & \small 1.03e-05  & \small  1.03e-05    & \small 1.03e-05  & \small 1.4e-09  & \small 6.0e-10\\

\small $0.3$              & \small  1.02e-05   & \small 1.02e-05  & \small  1.02e-05    & \small 1.02e-05  & \small 1.4e-09  & \small 6.0e-10 \\

\small $0.4$              & \small  9.93e-06   & \small 9.93e-06  & \small  9.93e-06    & \small 9.93e-06  & \small 1.5e-09  & \small 8.0e-10 \\

\small $0.5$              & \small  9.62e-06   & \small 9.62e-06  & \small  9.62e-06    & \small 9.62e-06  & \small 2.6e-09  & \small 1.8e-09 \\

\small $0.6$              & \small  2.73e-06   & \small 6.07e-06  & \small  6.93e-06    & \small 9.25e-06  & \small 1.9e-09  & \small 1.2e-09\\

\small $0.7$              & \small  6.67e-07   & \small 3.65e-06  & \small  4.75e-06    & \small 8.75e-06  & \small 1.4e-09  & \small 7.0e-10\\

\small $0.8$              & \small  1.58e-06   & \small 2.02e-06  & \small  2.93e-06    & \small 7.88e-06  & \small 9.0e-10  & \small 3.0e-10\\

\small $0.9$              & \small  1.08e-06   & \small 8.76e-07  & \small  1.37e-06    & \small 5.78e-06  & \small 5.5e-10  & \small 1.1e-09\\
\small $1.0$              & \small  0          & \small 0  & \small  0    & \small 1.10e-10  & \small 2.74e-11   & \small 3.6e-11\\
\bottomrule
\end{tabular}\\[2pt]
\end{table}

\begin{table}[htb]
\caption{The theoretical estimate errors  $\mathrm{TE}_N$ and comparison of the maximal absolute errors $\mathrm{ME}_N$ of present method and of other methods for Example 2.}
\label{table:1}
\newcommand{\m}{\hphantom{$-$}}
\newcommand{\cc}[1]{\multicolumn{21}{c}{#1}}
\renewcommand{\tabcolsep}{0.54pc} 
\renewcommand{\arraystretch}{1.1} 
\begin{tabular}{@{}llllllllll}
\toprule
$N$    & \small $\mathrm{TE}_N$   & \small $\mathrm{ME}_N$ & \small $N$    & \small $\mathrm{TE}_N$ & \small $\mathrm{ME}_N$  & \small $N$   & \small in \cite{Caglar09} & \small $N$   & \small in \cite{Chawla88}\\
\hline
\small $10$  & \small 6.9957e-05 & \small 1.0488e-05  & \small $16$ & \small 2.2413e-07 & \small 3.5041e-08  & \small $20$ & \small 3.1607e-05 & \small $16$ & \small 2.52e-03 \\

\small $12$  & \small 1.0042e-05 & \small 1.5380e-06  & \small $18$ & \small 3.2730e-08 & \small 5.4593e-09   & \small $40$ & \small 7.8742e-06 & \small $32$ & \small 1.83e-04   \\

\small $14$ & \small 1.4795e-06 & \small 2.3036e-07 & \small $20$  & \small 6.6210e-09 & \small 8.4075e-10  & \small $60$  & \small 3.5011e-06 & \small $64$ & \small 1.28e-05  \\
\bottomrule
\end{tabular}\\[2pt]
\end{table}

\vskip 10pt
{\bf \noindent Example 3.} Consider the following nonlinear SBVP in the study of steady-state oxygen diffusion in a spherical cell \cite{Babolian13,Khuri10,Wazwaz11,Caglar09}:
\begin{equation}  \label{example3-1}
\hskip 15pt
\displaystyle u''(x)+\frac{\alpha}{x} u'(x)=\frac{\delta u(x)}{u(x)+\mu}, \hskip 5pt  \delta >0, \hskip 5pt \mu >0,
\end{equation}
subject to the boundary conditions
\begin{equation}   \label{example3-2}
\hskip 15pt
u'(0)=0, \hskip 5pt 5u(1)+u'(1)=5,
\end{equation}
where $\delta$ and $\mu$ are often taken as $0.76129$ and $0.03119$, respectively. We take the value of $\alpha$ as 1, 2 and 3.

The Adomian polynomials of nonlinear term $f(x,u)=\frac{\delta u(x)}{u(x)+\mu}$ in this problem are computes as
\[
\begin{array}{lll}
\hskip 25pt \displaystyle A_0=\frac{\delta }{U(0)+\mu}U(0),\\
\hskip 25pt \displaystyle A_1=\frac{\delta \mu }{(U(0)+\mu)^2}U(1) ,\\
\hskip 25pt \displaystyle A_2=\frac{\delta \mu } {(U(0)+\mu)^2}  U(2) - \frac{ \delta \mu } {(U(0)+\mu)^3}  U^2(1),\\
\hskip 25pt \displaystyle A_3=\frac{\delta \mu } {(U(0)+\mu)^2}  U(3) - \frac{2\delta \mu } {(U(0)+\mu)^3}U(1)U(2)+\frac{ \delta \mu } {(U(0)+\mu)^4}  U^3(1),\\
\hskip 45pt \vdots
\end{array}
\]

Proceeding as before, we compute the approximate solution $u_{12,2}(x)$ for $N=12$ and $\alpha=2$, and show a comparison of the numerical results compared to the other existing methods in Table 6, from which one can see that the results of our computations are in good agreement with those ones obtained by the SGM \cite{Babolian13}, the BSDM \cite{Khuri10}, the VIM \cite{Wazwaz11} and the BSM \cite{Caglar09}.

Moreover, since there is no exact solution of this problem, we instead investigate the absolute residual error functions and the maximal error remainder parameters, which are the measures of how well the numerical solution satisfies the original problem (\ref{example3-1})-(\ref{example3-2}). The absolute residual error functions are
\[
\hskip 25pt \left | \mathrm{ER}_{N,\alpha}(x)\right |=\left | u_{N,\alpha}''(x)+\frac{\alpha}{x} u_{N,\alpha}'(x)-\frac{\delta u_{N,\alpha}(x)}{\mu+u_{N,\alpha}(x)}\right |, \hskip 5pt 0< x \leq 1,
\]
and the maximal error remainder parameters are
\[
\hskip 25pt \mathrm{MER}_{N,\alpha} = \max \limits_{0< x \leq 1} \left | \mathrm{ER}_{N,\alpha}(x)\right |.
\]
In Fig. 1, we plot the absolute residual error functions $|\mathrm{ER}_{N,2}(x)|$ for $N=2$ through 12 by step 2. Besides, the maximal error remainder parameters $\mathrm{MER}_{N,\alpha}$ for the same $N$ and $\alpha=1,2,3$ are listed in Table 7, from which it is interesting to point out that for a given $N$ the accuracy of our approximate solutions increases with the increase of $\alpha$. Moreover, Fig. 1 and Table 7 show clearly that the accuracy of our method is getting better as the approximation level is increasing for a fixed $\alpha$.
The logarithm plots of the value of $\mathrm{MER}_{2,\alpha}$ through $\mathrm{MER}_{12,\alpha}$ for $\alpha=1,2,3$ are displayed in Fig. 2, which demonstrates an approximately exponential rate of convergence for the obtained truncated series solutions and thus the presented method converges rapidly to the exact solution.

\begin{table}[htb]
\caption{Comparison of the approximate solutions for Example 3.}
\label{table:1}
\newcommand{\m}{\hphantom{$-$}}
\newcommand{\cc}[1]{\multicolumn{21}{c}{#1}}
\renewcommand{\tabcolsep}{0.83pc} 
\renewcommand{\arraystretch}{1.1} 
\begin{tabular}{@{}llllll}
\toprule
$x$                       & \small BSDM \cite{Khuri10} & \small BSM \cite{Caglar09}  & \small VIM \cite{Wazwaz11} & \small SGM \cite{Babolian13}   & \small Present method  \\
\hline
\small $0.0$              & \small 0.8284832948     & \small 0.8284832729   & \small 0.8284832761     & \small 0.8284832912  & \small 0.8284832870  \\
\small $0.1$              & \small 0.8297060968     & \small 0.8297060752   & \small 0.8297060781     & \small 0.8297060933  & \small 0.8297060890  \\
\small $0.2$              & \small 0.8333747380     & \small 0.8333747169   & \small 0.8333747193     & \small 0.8333747345  & \small 0.8333747303   \\
\small $0.3$              & \small 0.8394899183     & \small 0.8394898981   & \small 0.8394898996     & \small 0.8394899148  & \small 0.8394899106   \\
\small $0.4$              & \small 0.8480527887     & \small 0.8480527703   & \small 0.8480527701     & \small 0.8480527859  & \small 0.8480527816    \\
\small $0.5$              & \small 0.8590649275     & \small 0.8590649139   & \small 0.8590649108     & \small 0.8590649281  & \small 0.8590649239   \\
\small $0.6$              & \small 0.8725283156     & \small 0.8725283084   & \small 0.8725282997     & \small 0.8725283208  & \small 0.8725283166  \\
\small $0.7$              & \small 0.8884452994     & \small 0.8884452958   & \small 0.8884452781     & \small 0.8884453065  & \small 0.8884453023 \\
\small $0.8$              & \small 0.9068185417     & \small 0.9068185402   & \small 0.9068185095     & \small 0.9068185490  & \small 0.9068185448  \\
\small $0.9$              & \small 0.9276509830     & \small 0.9276509825   & \small 0.9276509392     & \small 0.9276509893  & \small 0.9276509853  \\
\small $1.0$              & \small 0.9509457948     & \small 0.9509457946   & \small 0.9509457539     & \small 0.9509457994  & \small 0.9509457960  \\
\bottomrule
\end{tabular}\\[2pt]
\end{table}

\begin{table}[htb]
\caption{The maximal error remainder parameters $\mathrm{MER}_{N,\alpha}$ for Example 3.}
\label{table:1}
\newcommand{\m}{\hphantom{$-$}}
\newcommand{\cc}[1]{\multicolumn{21}{c}{#1}}
\renewcommand{\tabcolsep}{0.85pc} 
\renewcommand{\arraystretch}{1.2} 
\begin{tabular}{@{}lllllll}
\toprule
$\alpha$       & \small $\mathrm{MER}_{2,\alpha}$   & \small $\mathrm{MER}_{4,\alpha}$  & \small $\mathrm{MER}_{6,\alpha}$ & \small $\mathrm{MER}_{8,\alpha}$   & \small $\mathrm{MER}_{10,\alpha}$  & \small $\mathrm{MER}_{12,\alpha}$ \\
\hline
\small $1$              & \small  5.8000e-03    & \small 1.4000e-03   & \small 3.1751e-04     & \small 7.3547e-05  & \small 1.7000e-05 & \small 3.9243e-06 \\
\small $2$              & \small  3.4000e-03    & \small 4.8431e-04   & \small 6.7761e-05     & \small 9.4474e-06  & \small 1.3142e-06 & \small 1.8267e-07 \\
\small $3$              & \small  2.4000e-03    & \small 2.4481e-04   & \small 2.4485e-05     & \small 2.4388e-06  & \small 2.4240e-07 & \small 2.4065e-08 \\
\bottomrule
\end{tabular}\\[2pt]
\end{table}

\vskip 10pt
{\bf \noindent Example 4.} Consider the following nonlinear SBVP which arises in the study of the distribution of heat sources in the human head \cite{Pandey97,Rashidinia07,Kanth06,Babolian13,Khuri10,Singh14,Caglar09,Duggan86}:
\begin{equation}
\hskip 15pt
\displaystyle u''(x)+\frac{2}{x} u'(x)=-e^{-u(x)},
\end{equation}
subject to the boundary conditions
\begin{equation}
\hskip 15pt
u'(0)=0, \hskip 5pt au(1)+bu'(1)=0.
\end{equation}

We consider the following two cases:

\textit {Case one: } $a=b=1.$

\textit {Case two: } $ a=0.1, b=1.$

The Adomian polynomials of nonlinear term $f(x,u)=-e^{-u(x)}$ in this problem are computed as
\[
\begin{array}{lll}
\hskip 25pt \displaystyle A_0=-e^{-U(0)},\\
\hskip 25pt \displaystyle A_1=U(1)e^{-U(0)},\\
\hskip 25pt \displaystyle A_2=U(2)e^{-U(0)}-\frac{1}{2}U^2(1)e^{-U(0)},\\
\hskip 25pt \displaystyle A_3=U(3)e^{-U(0)}-U(1)U(2)e^{-U(0)}+\frac{1}{6}U^3(1)e^{-U(0)},\\
\hskip 45pt \vdots
\end{array}
\]

Again no exact solution exists for this equation, hence it was handled numerically. Table 8 describes the numerical results of the first case obtained by the proposed method at the order of approximation $N=12$ and the other existing methods, including the FDM \cite{Pandey97}, the non-polynomial cubic spline method (NPCSM) \cite{Rashidinia07}, the CSM \cite{Kanth06} and the SGM \cite{Babolian13}. Meanwhile, a comparison for the approximate solutions of the second case obtained by the present method with the same approximation level as the first case and the previous existing methods which include the CSM \cite{Kanth06}, the SGM \cite{Babolian13}, the BSDM \cite{Khuri10} and the BSM \cite{Caglar09} is presented in Table 9. One can seen from two Tables that our computations are in good line with the results obtained by the other approaches compared. In fact, at the approximation level for $N=12$, the maximal absolute error is found to be order of magnitude O($10^{-7}$) for the first case, and O($10^{-9}$) for the second case.

\begin{table}[htb]
\caption{Comparison of the numerical results for the first case of Example 4.}
\label{table:1}
\newcommand{\m}{\hphantom{$-$}}
\newcommand{\cc}[1]{\multicolumn{21}{c}{#1}}
\renewcommand{\tabcolsep}{0.70pc} 
\renewcommand{\arraystretch}{1.1} 
\begin{tabular}{@{}llllll}
\toprule
$x$                       & \small FDM \cite{Pandey97} & \small NPCSM \cite{Rashidinia07}  & \small CSM \cite{Kanth06} & \small SGM \cite{Babolian13}   & \small Present method  \\
\hline
\small $0.0$              & \small 0.3675169710     & \small 0.3675181074   & \small 0.3675179806     & \small 0.3675168124  & \small 0.3675167997  \\
\small $0.1$              & \small 0.3663623697     & \small 0.3663637561   & \small 0.3663634922     & \small 0.3663623265  & \small 0.3663623137  \\
\small $0.2$              & \small 0.3628941066     & \small 0.3628959378   & \small 0.3628952219     & \small 0.3628940634  & \small 0.3628940507   \\
\small $0.3$              & \small 0.3570975862     & \small 0.3570991429   & \small 0.3570986892     & \small 0.3570975430  & \small 0.3570975301   \\
\small $0.4$              & \small 0.3489484612     & \small 0.3489499903   & \small 0.3489495462     & \small 0.3489484178  & \small 0.3489484049    \\
\small $0.5$              & \small 0.3384121893     & \small 0.3384136581   & \small 0.3384132502     & \small 0.3384121459  & \small 0.3384121330   \\
\small $0.6$              & \small 0.3254435631     & \small 0.3254450019   & \small 0.3254445925     & \small 0.3254435196  & \small 0.3254435063  \\
\small $0.7$              & \small 0.3099860810     & \small 0.3099878567   & \small 0.3099870705     & \small 0.3099860373  & \small 0.3099860240 \\
\small $0.8$              & \small 0.2919711440     & \small 0.2919789654   & \small 0.2919720836     & \small 0.2919711001  & \small 0.2919710864  \\
\small $0.9$              & \small 0.2713170512     & \small 0.2713185637   & \small 0.2713179289     & \small 0.2713170072  & \small 0.2713169936  \\
\small $1.0$              & \small 0.2479277646     & \small 0.2479292837   & \small 0.2479285659     & \small 0.2479277203  & \small 0.2479277073  \\
\bottomrule
\end{tabular}\\[2pt]
\end{table}

\begin{table}[htb]
\caption{Comparison of the numerical results for the second case of Example 4.}
\label{table:1}
\newcommand{\m}{\hphantom{$-$}}
\newcommand{\cc}[1]{\multicolumn{21}{c}{#1}}
\renewcommand{\tabcolsep}{0.85pc} 
\renewcommand{\arraystretch}{1.1} 
\begin{tabular}{@{}llllll}
\toprule
$x$                       & \small CSM \cite{Kanth06} & \small BSM \cite{Caglar09}  & \small BIDM \cite{Khuri10} & \small SGM \cite{Babolian13}   & \small Present method  \\
\hline
\small $0.0$              & \small 1.147041084      & \small 1.147039937    & \small 1.147040795     & \small 1.147039016   & \small 1.147039019  \\
\small $0.1$              & \small 1.146511706      & \small 1.146510559    & \small 1.146511419     & \small 1.146509639   & \small 1.146509642  \\
\small $0.2$              & \small 1.144922563      & \small 1.144921418    & \small 1.144922282     & \small 1.144920499   & \small 1.144920502   \\
\small $0.3$              & \small 1.142270622      & \small 1.142269478    & \small 1.142270348     & \small 1.142268560   & \small 1.142268563   \\
\small $0.4$              & \small 1.138550801      & \small 1.138549661    & \small 1.138550539     & \small 1.138548745   & \small 1.138548748    \\
\small $0.5$              & \small 1.133755950      & \small 1.133754813    & \small 1.133755703     & \small 1.133753900   & \small 1.133753904   \\
\small $0.6$              & \small 1.127876795      & \small 1.127875663    & \small 1.127876562     & \small 1.127874754   & \small 1.127874756  \\
\small $0.7$              & \small 1.120901889      & \small 1.120900762    & \small 1.120901665     & \small 1.120899858   & \small 1.120899860 \\
\small $0.8$              & \small 1.112817535      & \small 1.112816416    & \small 1.112817317     & \small 1.112815517   & \small 1.112815520  \\
\small $0.9$              & \small 1.103607704      & \small 1.103606593    & \small 1.103607490     & \small 1.103605701   & \small 1.103605704  \\
\small $1.0$              & \small 1.093253927      & \small 1.093252826    & \small 1.093253716     & \small 1.093251942   & \small 1.093251944   \\
\bottomrule
\end{tabular}\\[2pt]
\end{table}

\vskip 10pt
{\bf \noindent Example 5.} Consider the following SBVP with nonlinear term different from the forms (4)-(7) which arises in the radial stress on a rotationally symmetric shallow membrane cap \cite {Singh14,Kanth10}:
\begin{equation}  \label{5-1}
\hskip 15pt
\displaystyle u''(x)+\frac{3}{x} u'(x)=\frac{1}{2}-\frac{1}{8u^2(x)},
\end{equation}
subject to the boundary conditions
\begin{equation}
\hskip 15pt
u'(0)=0, \hskip 5pt u(1)=1.
\end{equation}

The Adomian polynomials of nonlinear term $f(x,u)=\frac{1}{2}-\frac{1}{8u^2(x)}$ in this problem are computed as
\[
\begin{array}{lll}
\hskip 25pt \displaystyle A_0=\frac{1}{2}-\frac{1}{8U^2(0)},\\
\hskip 25pt \displaystyle A_1=\frac{1}{4} \frac{U(1)}{U^3(0)},\\
\hskip 25pt \displaystyle A_2=-\frac{3}{8} \frac{U^2(1)}{U^4(0)}+ \frac{1}{4} \frac{U(2)}{U^3(0)},\\
\hskip 25pt \displaystyle A_3=\frac{1}{2} \frac{U^3(1)}{U^5(0)}- \frac{3}{4} \frac{U(1)U(2)}{U^4(0)}+ \frac{1}{4} \frac{U(3)}{U^3(0)},\\
\hskip 45pt \vdots
\end{array}
\]

Like the previous problems 3 and 4, a closed-form solution to this equation can not be written down. So we instead investigate the absolute residual error functions and the maximal error remainder parameters to examine the accuracy and the reliability of our numerical results. Here, the absolute residual error functions are
\[
\hskip 25pt \left | \mathrm{ER}_N(x) \right | = \left | u_N''(x)+\frac{3}{x} u_N'(x)-\frac{1}{2}+\frac{1}{8u_N^2(x)} \right |, \hskip 5pt 0 < x \leq 1,
\]
and the maximal error remainder parameters are
\[
\hskip 25pt \mathrm{MER}_N = \max \limits_{0< x \leq 1} \left | \mathrm{ER}_N(x) \right |.
\]
In Fig. 3, we plot the absolute residual error functions $\left | \mathrm{ER}_N(x) \right |$ for $N=4$ through 14 by step 2. The logarithm plot for the maximal error remainder parameters $\mathrm{MER}_N$ for the same $N$ is shown in Fig. 4, which demonstrates an approximately exponential rate of convergence of the obtained truncated series solutions and thus the presented method converges rapidly to the exact solution. Even though there is no exact solution for this problem, the following 10th order approximation has an accuracy of O($10^{-8}$) and can be used for practical applications
\[
\begin{array}{lll}
\hskip 25pt u_{10}(x)=0.9541353070+(0.4533672772e-1)x^2+(0.5436871104e-3)x^4-\\
\hskip 72pt (0.1611538997e-4)x^6+(0.3997114810e-6)x^8-\\
\hskip 72pt (0.6144814593e-8)x^{10}.
\end{array}
\]

\section{Conclusion}
In this work, a reliable approach based on the IDTM is presented to handle the numerical solutions of a class of nonlinear SBVPs arising in various physical models. This scheme takes the form of a truncated series with easily computable coefficients via the Adomian polynomials of those nonlinearities in the given problem. With the proposed algorithm, there is no need of discretization of the variables, linearization or small perturbation. Numerical results show that the proposed method works well for the SBVPs (1)-(3) with a satisfying low error. Besides, it is obvious that evaluation of more components of the approximate solution will reasonably improve the accuracy of truncated series solution by using the proposed method. Comparisons of the results reveal that the present method is very effective and accurate. Moreover, we are convinced that the IDTM can be extended to solve the other type of functional equations involving nonlinear terms more easily as the Adomian polynomials are applicable for any analytic nonlinearity and can be generated quickly with the aid of the algorithm proposed by Duan.

It is necessary to point out that algebraic equation (23) is a nonlinear one, and we shall inevitably encounter the bad roots while solving it. The criterion to separate the good root from a swarm of bad ones is convergence because it represents the value of unknown function at the origin and will not change for the different $N$.

\section*{Acknowledgement}

This work was supported by the Scientific Research Fund of Zhejiang Provincial Education Department of China (No.Y201430940) and K.C. Wong Magna Fund in Ningbo University, and partially supported by the National Nature Science Foundation of China (No. 11226243).

\newpage
{\bfseries List of Figures}\vspace{12pt}

Figure 1. The absolute residual error functions $\left | {\mathrm {ER}_{N,2}(x)} \right|$ for $N=2,4,6$ (left) and $8,10,12$ (right) of Example 3.\\

Figure 2. The logarithmic plots for the maximal error remainder parameters $\mathrm{MER}_{N,\alpha}$ for $N=2$ through $12$ by step $2$ and $\alpha=1$ (up, left), $\alpha=2$ (up, right), $\alpha=3$ (down) of Example 3.\\

Figure 3. The absolute residual error functions $\left | {\mathrm {ER}_N(x)} \right|$ for $N=4,6,8$ (left) and $10,12,14$ (right) of Example 5.\\

Figure 4. The logarithmic plot for the maximal error remainder parameters $\mathrm{MER}_N$ for $N=2$ through $14$ by step $2$ of Example 5.\\

\newpage
\begin{figure}[!hbp]
\centering
\includegraphics[width=15cm,height=5.2cm]{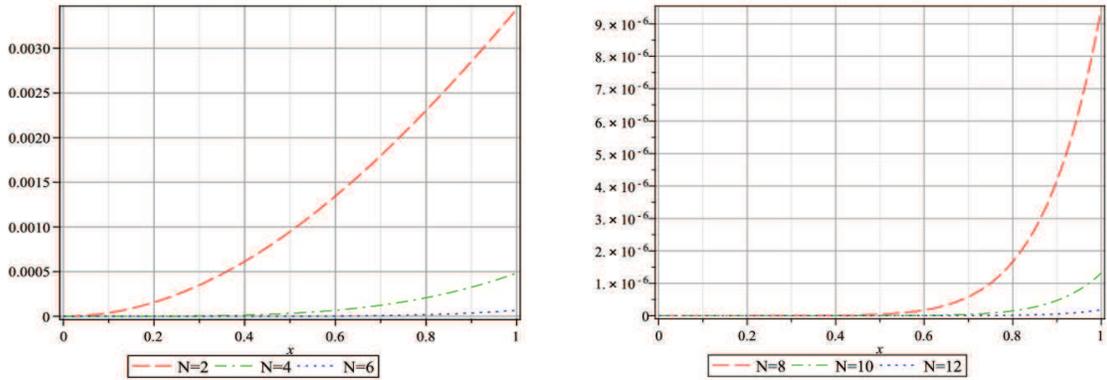}
\caption{The absolute residual error functions $\left | {\mathrm {ER}_{N,2}(x)} \right|$ for $N=2,4,6$ (left) and $8,10,12$ (right) of Example 3.}
\end{figure}

\begin{figure}[!hbp]
\centering
\includegraphics[width=15cm,height=10cm]{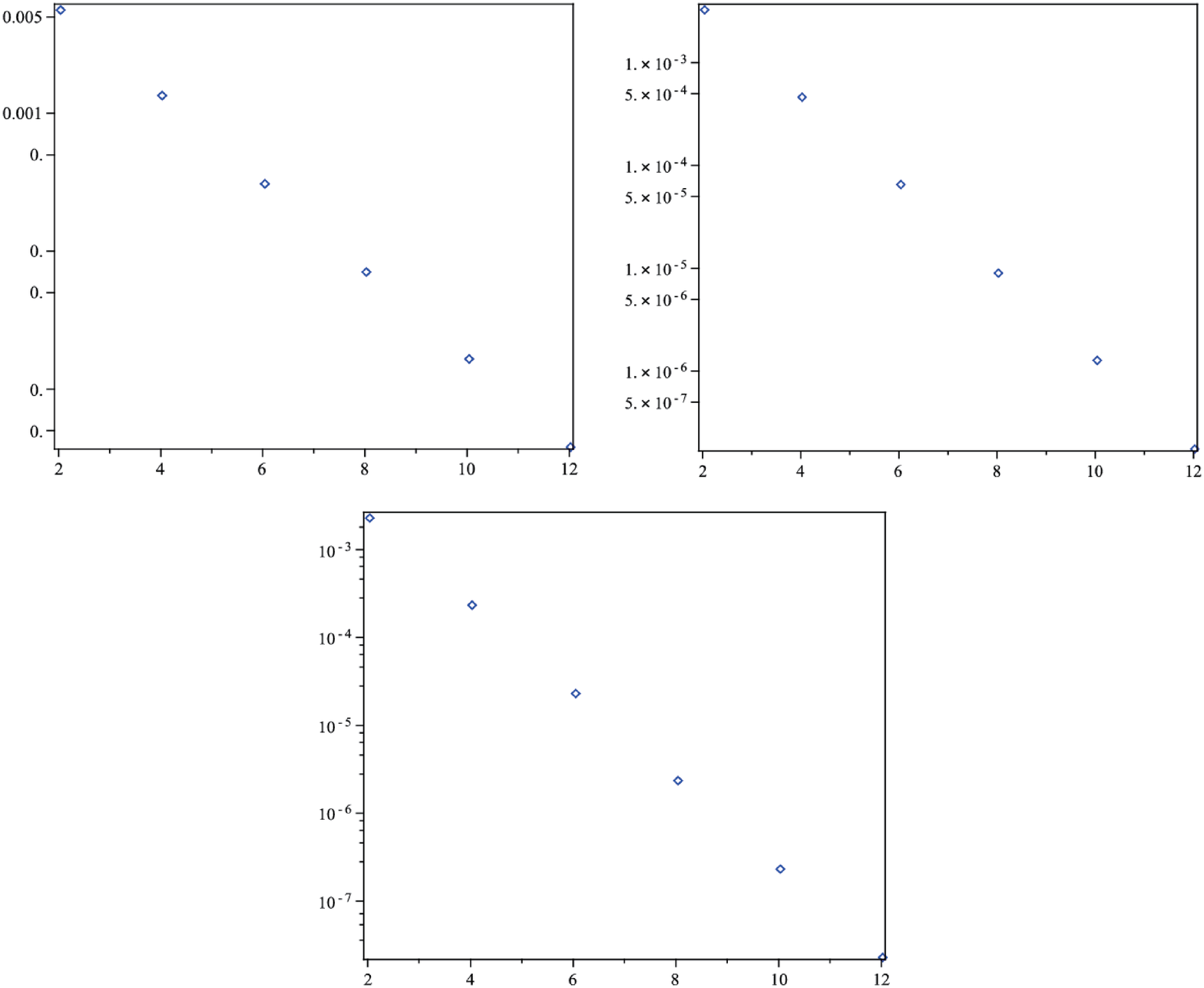}
\caption{The logarithmic plots for the maximal error remainder parameters $\mathrm{MER}_{N,\alpha}$ for $N=2$ through $12$ by step $2$ and $\alpha=1$ (up, left), $\alpha=2$ (up, right), $\alpha=3$ (down) of Example 3.}
\end{figure}

\begin{figure}[!hbp]
\centering
\includegraphics[width=16cm,height=6cm]{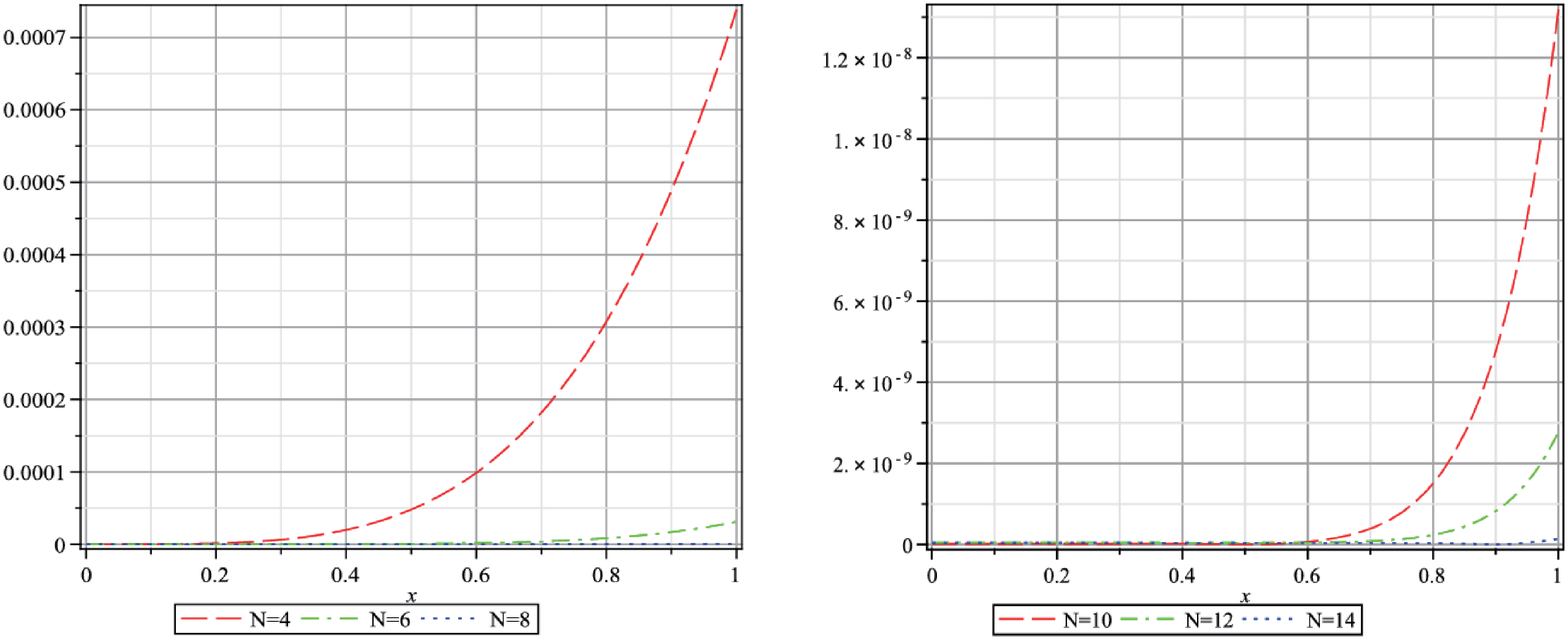}
\caption{The absolute residual error functions $\left | {\mathrm {ER}_N(x)} \right|$ for $N=4,6,8$ (left) and $10,12,14$ (right) of Example 5.}
\end{figure}

\begin{figure}[!hbp]
\centering
\includegraphics[width=10cm,height=8cm]{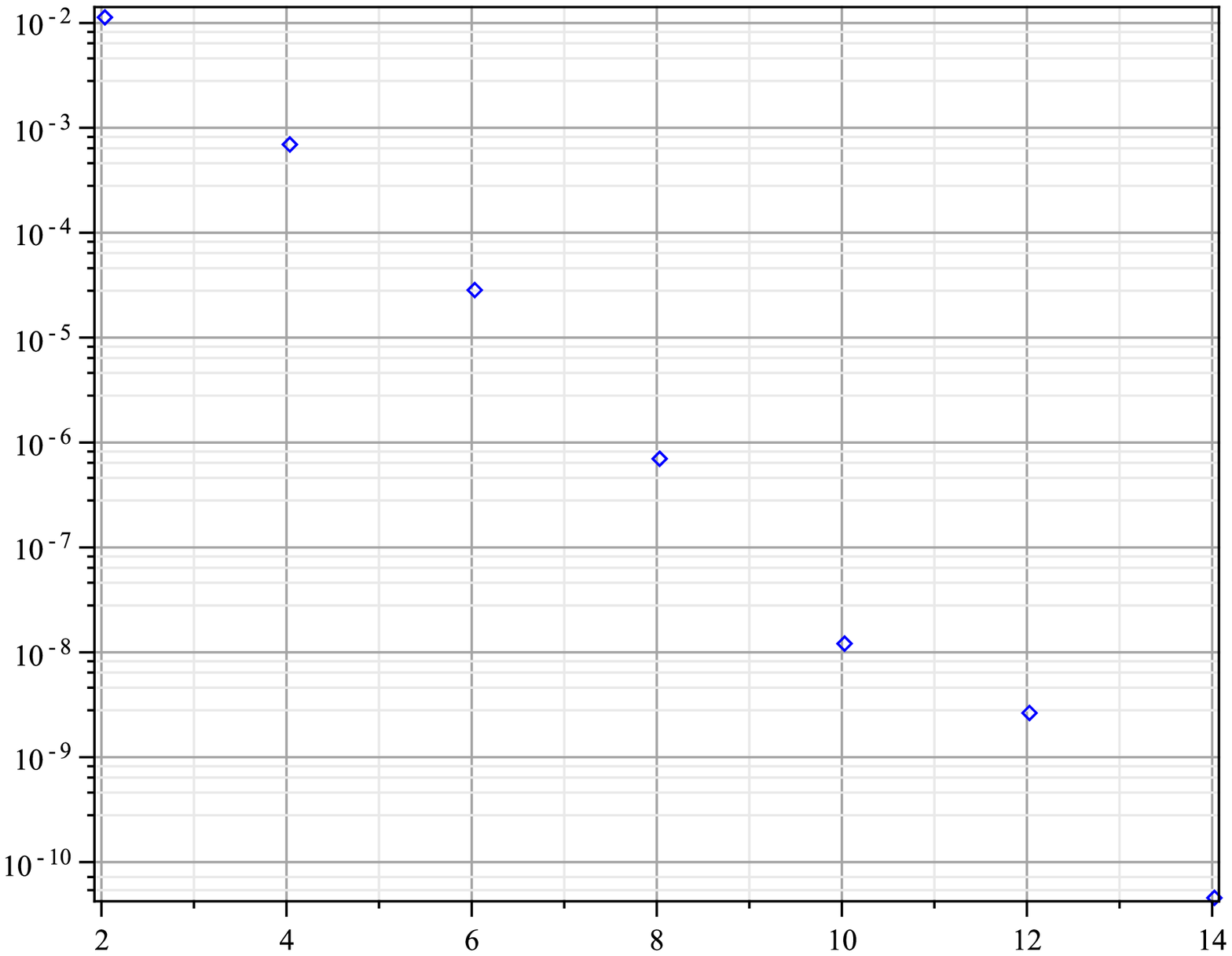}
\caption{The logarithmic plot for the maximal error remainder parameters $\mathrm{MER}_N$ for $N=2$ through $14$ by step $2$ of Example 5.}
\end{figure}


\begin{thebibliography}{99}

\bibitem{Russel75} R.D. Russell, L.F. Shampine, Numerical methods for singular boundary value problems, SIAM J. Numer. Anal. 12(1)(1975)13-36.

\bibitem{Lin76} S.-H. Lin, Oxygen diffusion in a spherical cell with nonlinear oxygen uptake kinetics, J. Theor. Biol. 60(2)(1976)449-457.
\bibitem{McElwain78} D.L.S. McElwain, A re-examination of oxygen diffusion in a spherical cell with michaelis-menten oxygen uptake kinetics, J. Theor. Biol. 71(2)(1978)255-263.
\bibitem{Hiltmann83} P. Hiltmann, P. Lory, On oxygen diffusion in a spherical cell with Michaelis-Menten oxygen uptake kinetics, Bull. Math. Biol. 45(5)(1983)661-664.
\bibitem{Anderson85} N. Anderson, A.M. Arthurs, Analytical bounding functions in a spherical cell with michaelis-menten oxygen uptake kinetics, Bull. Math. Biol. 47(1)(1985)145-153.
\bibitem{Pandey97} R.K. Pandey, A finite difference method for a class of singular two point boundary value problems arising in physiology, Int. J. Comput. Math. 65(1-2)(1997)131-140.
\bibitem{Rashidinia07} J. Rashidinia, R. Mohammadi, R. Jalilian, The numerical solution of non-linear singular boundary value problems arising in physiology, Appl. Math. Comput. 185(1)(2007)360-367.
\bibitem{Kanth06} A.S.V. Ravi Kanth, V. Bhattacharya, Cubic spline for a class of non-linear singular boundary value problems arising in physiology, Appl. Math. Comput. 174(1)(2006)768-774.
\bibitem{Babolian13} E. Babolian, A. Eftekhari, A. Saadatmandi, A Sinc-Galerkin technique for the numerical solution of a class of singular boundary value problems, Comp. Appl. Math. (2013)1-19.


\bibitem{Khuri10}  S.A. Khuri, A. Sayfy, A novel approach for the solution of a class of singular boundary value problems arising in physiology, Math. Comput. Modell. 52(3-4)(2010)626-636.
\bibitem{Wazwaz13} A.M. Wazwaz, R. Rach, J.-S. Duan, Adomian decomposition method for solving the Volterra integral form of the Lane-Emden equations with initial and boundary conditions, Appl. Math. Comput. 219(10)(2013)5004-5019.
\bibitem{Singh14} R. Singh, J. Kumar, An efficient numerical technique for the solution of nonlinear singular boundary value problems, Comput. Phys. Commun. 185(4)(2014)1282-1289.
\bibitem{Kanth10} A.S.V. Ravi Kanth, K. Aruna, He's variational iteration method for treating nonlinear singular boundary problems, Comput. Math. Appl. 60(3)(2010)821-829.
\bibitem{Wazwaz11} A.M. Wazwaz, The variational iteration method for solving nonlinear singular boundary value problems arising in various physical models, Commun. Nonlinear Sci. Numer. Simulat. 16(10)(2011)3881-3886.

\bibitem{Turkyilmazoglu13} M. Turkyilmazoglu, Effective computation of exact and analytic approximate solutions to singular nonlinear equations of Lane-Emden-Fowler type, Appl. Math. Modell. 37(14-15)(2013)7539-7548.
\bibitem{Caglar09} H. \c{C}a\u{g}lar, N. \c{C}a\u{g}lar, M. \"{O}zer, B-spline solution of non-linear singular boundary value problems arising in physiology, Chaos Soliton. Fract. 39(3)(2009)1232-1237.

\bibitem{Flesch75} U. Flesch, The distribution of heat sources in the human head: a theoretical consideration, J. Theor. Biol. 54(2)(1975)285-287.
\bibitem{Gray80} B.F. Gray, The distribution of heat sources in the human head --- theoretical considerations, J. Theor. Biol. 82(3)(1980)473-476.
\bibitem{Duggan86} R.C. Duggan, A.M. Goodman, Pointwise bounds for a nonlinear heat conduction model of the human head, Bull. Math. Biol. 48(2)(1986)229-236.

\bibitem{Kumar10} M. Kumar, N. Singh, Modified Adomian decomposition method and computer implementation for solving singular boundary value problems arising in various physical problems, Comput. Chem. Eng. 34(11)(2010)1750-1760.
\bibitem{Chang14}  S.-H. Chang, Taylor series method for solving a class of nonlinear singular boundary value problems arising in applied science, Appl. Math. Comput. 235(25)(2014)110-117.
\bibitem{Chang12} S.-H. Chang, Electroosmotic flow in a dissimilarly charged slit microchannel containing salt-free solution, Eur. J. Mech. B-Fluid. 34(2012)85-90.

\bibitem{Chawla88} M.M. Chawla, R. Subramanian, H.L. Sathi, A fourth order method for a singular two-point boundary value problem, BIT 28(1)(1988)88-97.

\bibitem{Chandrasekhar39} S. Chandrasekhar, An introduction to the study of stellar strcture, Dover, New York, 1939.
\bibitem{Ames68} W.F. Ames, Nonlinear ordinary differential equations in tranport process, Academic press, New York, 1968.


\bibitem{Pukhov80} G.E. Pukhov, Differential transforms of functions and equations, Naukova Dumka, Kiev, 1980 (in Russian).
\bibitem{Pukhov82} G.E. Pukhov, Differential transforms and circuit theory, Int. J. Circ. Theor. App. 10 (1982)265-276.
\bibitem{Pukhov86} G.E. Pukhov, Differential transformations and mathematical modeling of physical processes, Naukova Dumka, Kiev, 1986 (in Russian).
\bibitem{Zhou86} J.-K. Zhou, Differential transformation and its applications for electrical circuits, Wuhan China: Huazhong University Press, 1986 (in Chinese).

\bibitem{Xie15} L.-J. Xie, C.-L. Zhou, S. Xu, A new algorithm based on differential transform method for solving multi-point boundary value problems, Int. J. Comput. Math. (2015) 1-14. http://dx.doi.org/10.1080/00207160.2015.1012070.

\bibitem{Odibat10} Z.M. Odibat, C. Bertelle, M.A. Aziz-Alaoui, G.H.E. Duchamp, A multi-step differential transform method and application to non-chaotic or chaotic systems, Comput. Math. Appl. 59(4)(2010)1462-1472.
\bibitem{Gokdogan12} A. G\"{o}kdo\u{g}an, M. Merdan, A. Yildirim, The modified algorithm for the differential transform method to solution of Genesio systems, Commum. Nonlinear Sci. 17(1)(2012) 45-51.

\bibitem{Momani08} S. Momani, V.S. Ert\"{u}rk, Solutions of non-linear oscillators by the modified differential transform method, Comput. Math. Appl. 55(4)(2008)833-842.
\bibitem{Chang08} S.-H. Chang, I-L. Chang, A new algorithm for calculating one-dimensional differential transform of nonlinear functions, Appl. Math. Comput. 195(2)(2008)799-808.



\bibitem{Ebaid12} A. Elsaid, Fractional differential transform method combined with the Adomian polynomials, Appl. Math. Comput. 218(12)(2012)6899-6911.
\bibitem{Fat13} H. Fatoorehchi, H. Abolghasemi, Improving the differential transform method: A novel technique to obtain the differential transforms of nonlinearities by the Adomian polynomials, Appl. Math. Modell. 37(8)(2013)6008-6017.

\bibitem{Duan10-6} J.-S. Duan, Recurrence triangle for Adomian polynomials, Appl. Math. Comput. 216(4)(2010)1235-1241.
\bibitem{Duan10-7} J.-S. Duan, An efficient algorithm for the multivariable Adomian polynomials, Appl. Math. Comput. 217(6)(2010)2456-2467.
\bibitem{Duan11} J.-S. Duan, Convenient analytic recurrence algorithms for the Adomian polynomials, Appl. Math. Comput. 217(13)(2011)6337-6348.

\bibitem{Adomian90} G. Adomian, A review of the decomposition method and some recent results for nonlinear equations, Math. Comput. Modell. 13(7)(1990)17-43.
\bibitem{Adomian94} G. Adomian, Solving Frontier Problems of Physics: The Decomposition Method, Kluwer Academic: Dordrecht, 1994.



\bibitem{Adomian83} G. Adomian, R. Rach, Inversion of nonlinear stochastic operators, J. Math. Anal. Appl. 91(1)(1983)39-46.
\bibitem{Rach08} R. Rach, A new definition of the Adomian polynomials, Kybernetes 37(7)(2008)910-955.
\bibitem{Rach84} R. Rach, A convenient computational form for the Adomian polynomials, J. Math. Anal. Appl. 102(2)(1984)415-419.
\bibitem{Wazwaz00} A.M. Wazwaz, A new algorithm for calculating Adomian polynomials for nonlinear operators, Appl. Math. Comput. 111(1)(2000)33-51.
\bibitem{Abbaoui95} K. Abbaoui, Y. Cherruault, V. Seng, Practical formulae for the calculus of multivariable Adomian polynomials, Math. Comput. Modell. 22(1)(1995)89-93.
\bibitem{Abdelwahid03} F. Abdelwahid, A mathematical model of Adomian polynomials, Appl. Math. Comput. 141(2-3)(2003)447-453.
\bibitem{Azreg09} M. Azreg-A\"{I}nou, A developed new algorithm for evaluating Adomian polynomials, CMES-Comput. Model. Eng. 42(1)(2009)1-18.

\end{thebibliography}
\end{document}